\documentclass{icmart}
\contact[minicozz@jhu.edu]{Department of Mathematics\\
Johns Hopkins University\\
3400 N. Charles St.\\
Baltimore, MD 21218}

\newtheorem{theorem}{Theorem}[section]

\theoremstyle{definition}

\newtheorem{remark}[theorem]{Remark}

\def\RR{{\bold R}}
\def\SS{{\bold S}}

\newcommand{\nn}{{\bold {n}}}

\def\CC{{\bold C }}
\newcommand{\eqr}[1]{(\ref{#1})}

\newcommand{\ee}{{\operatorname {e}}}
\newcommand{\cB}{{\mathcal {B}}}
\newcommand{\cP}{{\mathcal {P}}}
\newcommand{\cS}{{\mathcal {S}}}
\newcommand{\cSt}{{\mathcal{S}_{neck}}}
\newcommand{\cSu}{{\mathcal{S}_{ulsc}}}

\newcommand{\cI}{{\mathcal {I}}}
\newcommand{\cL}{{\mathcal {L}}}

\newcommand{\cM}{{\mathcal {M}}}
\newcommand{\cF}{{\mathcal {F}}}
\newcommand{\cT}{{\mathcal {T}}}
\newcommand{\barga}{\Gamma_{Clos}}

\title[Embedded minimal surfaces]{Embedded minimal surfaces}

\author[William P. Minicozzi II]{William P. Minicozzi II\thanks{The author was partially supported by NSF Grant DMS
 0405695.}}

\begin{document}

\begin{abstract}
  The study of
embedded minimal surfaces in $\RR^3$ is a classical problem,
dating to the mid 1700's, and many people have made key
contributions. We will survey  a few recent advances, focusing on
joint work with Tobias H. Colding of MIT and Courant, and taking
the opportunity to focus on results that have not been highlighted
elsewhere.
\end{abstract}



\maketitle

\section{Introduction} \label{s:s0}

An immersed surface $\Sigma$ in $\RR^3$ is said to be
{\emph{minimal}} if it has zero mean curvature and is
{\emph{embedded}} if the immersion is injective.   The study of
embedded minimal surfaces in $\RR^3$ is a classical problem,
dating to the mid 1700's, and many people have made key
contributions.

Many of the recent results have been surveyed elsewhere and we
will take the opportunity to highlight results that have not been
as well covered, concentrating on recent joint work with Tobias H.
Colding of MIT and the Courant Institute.  We will also briefly
cover recent important results of W. Meeks and H. Rosenberg and of
W. Meeks, J. Perez, and A. Ros.    We refer to the following
  surveys for other perspectives:
\begin{itemize}
\item  For more on the structure of properly embedded minimal
surfaces, see the joint expository article \cite{CM10} with Tobias
H. Colding as well as the surveys \cite{MeP} of W. Meeks and J.
Perez, \cite{Pz} of J. Perez, \cite{Ro2} of H. Rosenberg, as well
as the joint surveys \cite{CM15}, \cite{CM18}, and \cite{CM19}
with Tobias H. Colding. \item For the construction of embedded
minimal surfaces, see the surveys \cite{HoWeWo2} of D. Hoffman, M.
Weber, and M. Wolf, \cite{Ka} of N. Kapouleas, and \cite{Tr} of M.
Traizet. \item For properness of minimal surfaces and the
Calabi-Yau Conjectures, see the paper \cite{CM16} as well as the
surveys \cite{Ma} of F. Martin, \cite{CM18}, \cite{CM19}, and
\cite{Pz}.
\end{itemize}

\subsection{Embedded minimal surfaces of fixed genus}
We have chosen to  concentrate  on the following central question:
\begin{itemize}

\item Can one compactify the space of embedded minimal
surfaces of fixed genus?

\end{itemize}
Roughly speaking, we show in \cite{CM7} that a sequence of
embedded minimal surfaces with fixed genus has a subsequence that
converges away from a singular set to a collection of parallel
planes. The precise structure of the singular set and of the
surfaces near the singular set depends on the topology of the
surfaces. Consequently,  we consider three separate cases:
\begin{enumerate}
\item When the surfaces are disks.
\item When the surfaces are (non-simply connected) planar domains; i.e., the case of genus zero.
\item When the surfaces have a fixed non-zero genus.
\end{enumerate}
The case of disks was completed in \cite{CM3}, \cite{CM4},
\cite{CM5} and \cite{CM6} and plays a key role in the other two
cases as well; the case of disks was surveyed in \cite{CM15} and
\cite{CM18}. The other two cases, which were completed in
\cite{CM7}, will be one of the focal points of this survey.

 A key step in the compactness results for embedded
minimal surfaces of fixed genus is a structure result that
describes what these surfaces look like. We have chosen to focus
on the compactness theorems rather than the underlying structure
results, largely because it serves as a unifying theme and allows
us to simplify some of the statements. Roughly speaking, two main
structure theorems for (non-simply connected) embedded minimal
planar domains
 from \cite{CM7} are:
\begin{itemize}
\item Any such surface {\underline{without}} small necks can be
obtained by gluing together two oppositely--oriented double spiral
staircases. \item Any such surface  {\underline{with}} small necks
can be decomposed into ``pairs of pants'' by  cutting the surface
along a collection of short curves. After the cutting, we are left
with graphical pieces that are defined over a disk with either one
or two sub--disks removed (a topological disk with two sub--disks
removed is called a {\emph{pair of pants}}).
\end{itemize}
Both of these structures occur as different extremes in the
two-parameter family of minimal surfaces known as the Riemann
examples.

\subsection{Embedded minimal annuli}
The simplest example of a non-simply connected planar domain is of
course an annulus.   In  \cite{CM9}, we obtained a precise
description of what an embedded minimal annulus in a ball must
look like -   roughly speaking, it must look like catenoid. This
illustrates a few of the ideas for the general pair of pants
decomposition of \cite{CM7} in a relatively simple setting.
  This description can be thought of as an
effective version of the main theorems of \cite{Co} and
\cite{CM14}; i.e., \cite{CM9} applies to an annulus $\Sigma$ with
$\partial \Sigma \subset \partial B_r(0)$   and as $r$ goes to
infinity we recover the results of \cite{Co} and \cite{CM14}.

\subsection{Properness and removable singularities}
The next result that we will highlight is the proof of
``properness'' in \cite{CM8}.  This properness was used in
\cite{CM6} to analyze a neighborhood of each singular point,
showing that an entire neighborhood is foliated by limit planes.
This can be viewed as a removable singularity theorem for minimal
laminations.  The proof of properness in \cite{CM8} works only in
the global case where we have a sequence of embedded minimal disks
in a sequence of expanding balls whose radii tend to infinity -
the local case is where the disks are in a fixed ball.  Perhaps
surprisingly, it turned out that properness can fail in the local
case: In the local case, we can get limits with non-removable
singularities. One local example with non-removable singularities
is constructed in \cite{CM11}.

\subsection{The global structure of
complete embedded minimal surfaces in $\RR^3$}
 As mentioned above, there have been many important
recent developments in the field.  We will survey two of these
where the results of \cite{CM3}--\cite{CM7}  play a role:
\begin{itemize}
\item The uniqueness of the helicoid, proven by W. Meeks and H.
Rosenberg in \cite{MeRo1}. \item  The curvature bound for embedded
minimal planar domains with bounded horizontal flux proven by W.
Meeks, J. Perez, and A. Ros in \cite{MePRs1}.
\end{itemize}
The uniqueness of the helicoid solved a long-standing problem that
was largely considered unapproachable until recently and also has
many applications.  We will sketch the proof and explain how the
lamination theorem and one-sided curvature estimate played a key
role.

The curvature bound of \cite{MePRs1} was the key step in solving
an old conjecture of J. Nittsche and an important step for
understanding the moduli space of embedded minimal planar domains.
We will explain the result, give an idea why it should be true,
and explain how the compactness theorems of \cite{CM7} play a role
in the proof.

We should point out that there is a key distinction between these
two results and the other results that we have discussed: these
results both use in an essential way that the surfaces are
complete and without boundary.

\section{Minimal surfaces}
An immersed surface $\Sigma \subset \RR^3$ is {\it{minimal}} if it
is a critical point for area, i.e., if it has zero mean curvature.
The {\emph{mean curvature}} is the trace of the second fundamental
form $A$; recall that the eigenvalues of $A$ are called the
{\emph{principal curvatures}}.  Our surface $\Sigma$ will always
be embedded and will have a well-defined unit normal $\nn$.  The
map
\begin{equation}
\nn : \Sigma \to \SS^2
\end{equation}
 is called the {\emph{Gauss map}}.
Note that $A$ is the differential of the Gauss map.

Observe that if $\Sigma \subset \RR^3$ is minimal, then so is
every rigid motion of $\Sigma$.  Furthermore, so is a dilation of
$\Sigma$, i.e., so is the surface
\begin{equation}
    \lambda \, \Sigma = \{ \lambda \, x \, | \, x \in \Sigma \} \,
    .
\end{equation}
This is because dilating $\Sigma$ by $\lambda$ dilates the second
fundamental form by $\lambda^{-1}$.

Note that minimal surfaces are not necessarily area-minimizing. A
surface is   {\it{stable}} if it satisfies the second derivative
test; obviously, area-minimizing surfaces are stable.

\subsection{Classical minimal surfaces}

The simplest example of a minimal surface is   a flat plane (where
the unit normal is constant and, hence, where $A = 0 $).

The only non-trivial rotationally invariant minimal surface is the
{\emph{ catenoid}} (discovered in 1776), i.e.,  the minimal
surface in $\RR^3$ parametrized by
\begin{equation}    \label{e:cat}
 (\cosh s\, \cos t,\cosh s\, \sin t,s)  \text{ where }s,\,t\in \RR\,
 .
\end{equation}
More precisely, since dilations preserve minimality, there is
 a one-parameter family of catenoids (modulo rigid motions) given by
\begin{equation}    \label{e:cat2}
 \lambda \, (\cosh s\, \cos t,\cosh s\, \sin t,s)  \text{ where }s,\,t\in \RR\,
 .
\end{equation}

The {\emph{ helicoid}} (also discovered in 1776) is the minimal
surface $\Sigma$ in $\RR^3$ parametrized by
\begin{equation}
(s\cos t,s\sin t,t)\text{ where }s,\,t\in \RR\, .
\end{equation}
Note that the helicoid is a ``double-spiral staircase'',
consisting of a straight line in each horizontal plane where these
lines rotate at constant speed.  It can also be thought of as the
union of the ``graphs''  of the functions $\theta$ and $\theta +
\pi$ together with the vertical axis.  We will make this last
characterization more precise later when we introduce the notion
of a multi-valued graph.

 The catenoid can be thought of as ``two planes glued together
along a small neck.''  Surprisingly, by a theorem of F. Lopez and
A. Ros, it is impossible to glue together any other
{\emph{finite}} number of planes to get a complete properly
embedded minimal planar domain.  However,  the {\emph{Riemann
examples}} (constructed by Riemann around 1860) give a periodic
collection of horizontal planes glued together along small necks.
This is actually (modulo rigid motions) a two parameter family of
surfaces, where the parameters can roughly be thought of as
\begin{itemize} \item the size of the necks (or injectivity
radius), and \item the angle from one to the next.
\end{itemize}
As the angle goes to zero, the necks get further and further apart
and the family degenerates to a collection of catenoids.  As the
angle goes to $\pi/2$, the necks become virtually on top of each
other and the family degenerates to the union of two oppositely
oriented helicoids.  There are very pretty pictures of this
available from David Hoffman's web page:

\vskip1mm \noindent
http://www.msri.org/about/sgp/jim/geom/minimal/library/riemann/index.html

\section{Embedded minimal surfaces with fixed genus}

As mentioned, we will focus on compactness theorems for a sequence
$\Sigma_i \subset \RR^3$ of embedded minimal surfaces. There are
various notions of weak convergence (e.g., as currents or
varifolds). However, for us, the sequence $\Sigma_i$ converges to
a surface $\Sigma_{\infty}$ at a point $x \in \Sigma$ if there is
a ball $B_{r}(x)$ so that:
\begin{itemize}
\item
 For every $i$ sufficiently large,  $B_r (x) \cap \Sigma_i$ is a
 (connected) graph over a (subset of) the tangent plane
 $T_x\Sigma_{\infty}$ of a function $u_i$.
 \item
 As $i \to \infty$, the functions $u_i$ converge smoothly to a
 function $u_{\infty}$ where $B_r (x) \cap \Sigma_{\infty}$ is the
 graph of $u_{\infty}$.
\end{itemize}
Notice that there are two obvious necessary conditions for the
sequence $\Sigma_i$ to converge in this sense: The curvatures and
areas of the sequence must be locally bounded.

 It is not hard to see that the lack of a local
area bound is not such a serious problem as long as we have
embeddedness. Namely, if we have a uniform curvature bound near
$x$, then  the components of $B_r (x) \cap \Sigma_i$ are
well-approximated by their tangent planes for $r$ small.
Embeddedness then implies that all of these tangent planes must be
almost parallel. In particular, these components are all graphs
over the same plane of functions with a uniform $C^1$ bound. We
can use the Arzela-Ascoli theorem to pass to a subsequence that
``converges'' to a {\underline{collection}} of minimal surfaces
that do not cross. The strong maximum principle then implies that
two of these limit surfaces must be identical if they touch at
all, i.e., they are like the leaves of a foliation. This sort of
structure is called a {\emph{lamination}}.

The failure of the curvature bound is a more serious problem and
will force us to allow for a singular set where the sequence
simply does not converge smoothly.  The simplest example of this
is a sequence  of rescalings $\lambda_i \, \Sigma$ with $\lambda_i
\to 0$ of a fixed non-flat complete embedded minimal surface
$\Sigma$. This scales the curvature by the factor $\lambda_i^{-1}$
and, thus, will force the curvature to blow up at the origin. For
example, a sequence of rescaled  catenoids converges with
multiplicity two to the punctured plane.  The convergence is
smooth except at $0$ where $|A|^2 \to \infty$.  Notice that $0$ is
a removable singularity for the limit.

It follows from Choi and Schoen, \cite{CiSc}, that a similar
singular compactness result holds as long as we assume a uniform
bound on the total curvature:
\begin{quote}
A subsequence converges smoothly with finite multiplicity away
from a finite set of singular points; these singular points are
then removable singularities for the limit surface.{\footnote{In
fact, one can say a good deal more about the convergence and the
structure of the limit; see the 1995 paper of A.~Ros in Indiana
Math.}}
\end{quote}

The situation is more complicated when there is no a priori total
curvature bound. For example, if we take a sequence  of rescaled
helicoids, then the curvature blows up along the entire vertical
axis but is bounded away from this axis.  Thus, we get that
\begin{itemize} \item The intersection of the rescaled helicoids
with a ball {\underline{away from}} the vertical axis gives a
collection of graphs over the plane $\{ x_3 = 0 \}$. As $i \to
\infty$, these graphs become flat and horizontal.
\item The intersection of the rescaled helicoids with a ball
{\underline{centered on}} the vertical axis gives a double spiral
staircase, rotating faster and faster as $i \to \infty$.
\end{itemize}
In particular, the sequence of rescaled helicoids converges away
from the vertical axis to a foliation by flat parallel planes.

\begin{remark}
The same thing happens when one rescales any surface asymptotic to
the helicoid - such as the genus one helicoid constructed by D.
Hoffman, M. Weber, and M. Wolf in \cite{HoWeWo1}.
\end{remark}

If we do the same rescaling to a fixed surface in the family of
Riemann examples, then we get convergence away from a line to a
foliation by horizontal planes. In this case, the line is
{\underline{not}} perpendicular to the planes.

However, unlike the catenoid and helicoid, the Riemann examples
are a two-parameter family.  By choosing the two parameters
appropriately, one can produce sequences of Riemann examples that
illustrate both of the two structure theorems:
\begin{enumerate}
\item
If we take a sequence of Riemann examples where the neck size is
fixed and the angles go to $\frac{\pi}{2}$, then the surfaces with
angle near $\frac{\pi}{2}$ can be obtained by gluing together two
oppositely--oriented double spiral staircases. Each double spiral
staircase looks like a helicoid. This sequence of Riemann examples
converges to a foliation by parallel planes.  The convergence is
smooth away from the axes of the two helicoids (these two axes are
the singular set   where the curvature blows up).
\item
Suppose now that we take a sequence of examples where the neck
sizes go to zero. In this case, the surfaces can be cut along
short curves into collections of graphical pairs of pants.  The
short curves converge to singular points where the curvature blows
up and the graphical pieces converge to flat planes except at
these points.
\end{enumerate}

\section{\cite{CM7}: Compactness of embedded minimal surfaces with fixed genus}

We turn next to the main compactness results of \cite{CM7} for
embedded minimal surfaces with fixed genus.   We will restrict our
discussion to the case of planar domains, i.e., when the surfaces
have genus zero, to simplify things. In any case, the general case
of fixed genus requires only minor changes.

{\bf{In this section, $\Sigma_i \subset B_{R_i} \subset \RR^3$ is
a sequence of compact embedded minimal planar domains with
$\partial \Sigma_i \subset
\partial B_{R_i}$. Moreover, we will assume that $R_i \to
\infty$.}}

\vskip2mm
 The singular set $\cS$ is defined to be the set of
points where the curvature is blowing up.  That is, a point $y$ in
$\RR^3$ is in $\cS$ for a sequence $\Sigma_i$ if
\begin{equation}
    \sup_{B_{r}(y)\cap
\Sigma_i}|A|^2\to\infty {\text{ as $i \to \infty$ for all $r>0$}}
.
\end{equation}
It is not hard to see that we can pass to a subsequence so that
$\cS$ is well-defined and, furthermore, if $x \notin \cS$, then
there exists $r_x > 0$ so that
\begin{equation}    \label{e:cbound}
    \sup_{i} \, \sup_{B_{r_x}(x) \cap \Sigma_i} \, |A| < \infty \,
    .
\end{equation}

\subsection{The finer structure of $\cS$: Where the topology concentrates}
Sequences of planar domains which are not simply connected are,
after passing to a subsequence, naturally divided into two
separate cases depending on whether or not the topology is
concentrating at points.  To distinguish between these cases, we
will say that a sequence of surfaces $\Sigma_i^2\subset \RR^3$ is
{\it{uniformly locally simply connected}} (or ULSC)  if for each
  $x \in \RR^3$, there exists a constant $r_0 > 0$
(depending on $x$) so that for every surface $\Sigma_i$
\begin{equation}    \label{e:ulsc2}
 {\text{each connected component of }} B_{r_0}(x) \cap \Sigma_i {\text{ is
 a disk.}}
\end{equation}
  For instance, a sequence of rescaled catenoids
   where the necks shrink to zero is not ULSC, whereas a
sequence of rescaled helicoids is.

Another way of locally distinguishing sequences where the topology
does not concentrate from sequences where it does comes from
analyzing the singular set.  The singular set $\cS$ consists of
two types of points. The first type is roughly modelled on
rescaled helicoids
 and the second on
rescaled catenoids:
\begin{itemize}
\item
A point $y$ in $\RR^3$ is in $\cSu$ if the curvature for the
sequence $\Sigma_i$ blows up at $y$ and the sequence is ULSC in a
neighborhood of $y$.
\item
A point $y$ in $\RR^3$ is in $\cSt$ if the sequence is not ULSC in
any neighborhood of $y$. In this case, a sequence of closed
non-contractible curves $\gamma_i \subset \Sigma_i$ converges to
$y$.
\end{itemize}
The sets $\cSt$ and $\cSu$ are obviously disjoint and
  the curvature blows up at both, so   $\cSt \cup \cSu \subset \cS$.   An easy
argument  (proposition $I.0.19$ in \cite{CM5}) shows that, after
passing to a further subsequence, we can assume that
\begin{equation}    \label{e:}
    \cS=  \cSt \cup \cSu \, .
\end{equation}
Note that $\cSt = \emptyset$ is equivalent to that the sequence is
ULSC as is the case for sequences of rescaled helicoids.  On the
other hand, $\cSu = \emptyset$ for sequences of rescaled
catenoids.  (These definitions of $\cSu$ and $\cSt$ are specific
to the genus zero case that we are focusing on now; the
definitions in the fixed genus case can be found in section $1.1$
of \cite{CM7}.)

\subsection{Compactness away from $\cS$}

If we combine the local curvature bound \eqr{e:cbound} away from
$\cS$ and a variation on the Arzela-Ascoli theorem, we can pass to
a subsequence so that the $\Sigma_i$'s converge away from $\cS$ to
a limit lamination $\cL'$ of $\RR^3 \setminus \cS$.

The leaves of $\cL'$ are smooth, but not necessarily complete,
surfaces. To make this precise, we define the closure $\barga$ of
a leaf $\Gamma$ of $\cL'$ to be the union of the closures of all
bounded (intrinsic) geodesic balls in $\Gamma$; that is, we fix a
point $x_{\Gamma} \in \Gamma$ and set
\begin{equation}    \label{e:closureaa}
    \barga = \bigcup_{r} \overline{ \cB_r (x_{\Gamma}) } \,
    ,
\end{equation}
where $\overline{ \cB_r (x_{\Gamma}) }$ is the closure of $\cB_r
(x_{\Gamma})$ as a subset of $\RR^3$.

Clearly, a leaf $\Gamma$ is complete if and only if $\barga =
\Gamma$ and we always have that
\begin{equation}
    \barga \setminus \Gamma \subset \cS \, .
\end{equation}
The incomplete leaves of $\Gamma$ can be divided into several
types, depending on how $\barga$ intersects $\cS$:
\begin{itemize}
\item Collapsed leaves  where $\barga \cap \cSu$ contains a
removable singularity for $\Gamma$. \item Leaves $\Gamma$ with
$\barga \cap \cSu \ne \emptyset$, but where $\Gamma$ does not have
a removable singularity.  This would occur, for example, if
$\Gamma$ spirals infinitely into the collapsed leaf through
$\barga \cap \cSu$. (We  show in \cite{CM7} that this does not
occur.)
\item Leaves $\Gamma$ where $\barga \setminus \Gamma \subset
\cSt$; these obviously don't occur in the ULSC case.
\end{itemize}

\subsection{Disks}

 Before discussing the general ULSC case, it is useful to
recall the case of disks.  One consequence of
\cite{CM3}--\cite{CM6} is that there are only two
{\underline{local models}} for ULSC sequences of embedded minimal
surfaces.  That is, locally in a ball in $\RR^3$, one of following
holds:
\begin{itemize}
\item  The curvatures are bounded and the surfaces are locally
{\underline{graphs}} over a plane. \item The curvatures blow up
and the surfaces are locally {\underline{double spiral
staircases}}.
\end{itemize}
Both of these cases are illustrated by taking a sequence of
rescalings of the helicoid; the first case occurs away from the
axis, while the second case occurs on the axis.

Using in part this local description, we were able to prove that
 any sequence of embedded
minimal disks with curvatures blowing up has a subsequence that
converges to a foliation by parallel planes.  This convergence is
away from a Lipschitz curve $\cS$ that is transverse to the
planes. (See the appendix for the precise statements.)

\subsection{Planar domains: the general structure theorems}

We will show that every sequence $\Sigma_i$ has a subsequence that
is either    ULSC or for which $\cSu$ is empty. This is the
 next ``no mixing'' theorem.  We will see later that these two different cases give
 two very different structures.

\begin{theorem}[No mixing theorem, \cite{CM7}]     \label{c:main}
If the $\Sigma_i$'s are genus zero, then there is a subsequence
with either $\cSu = \emptyset$ or $\cSt = \emptyset$.
\end{theorem}

Common for both the ULSC case and the case where $\cSu$ is empty
is that the limits are always  laminations by flat parallel planes
and  the singular sets are always  closed subsets contained in the
union of the planes.  This is the content of the next theorem:

\begin{theorem}[Planar lamination theorem, \cite{CM7}] \label{t:tab}
 If the $\Sigma_i$'s are genus zero and
\begin{equation}
\sup_{B_1\cap \Sigma_i}|A|^2\to \infty \, ,
 \end{equation}
   then there exists a subsequence $\Sigma_j$,
     a lamination $\cL=\{x_3=t\}_{ \{ t \in \cI \} }$ of $\RR^3$
 by parallel
 planes (where $\cI \subset \RR$ is a closed set), and a closed nonempty set
 $\cS$ in the union of the leaves of $\cL$ such that
 after a rotation of $\RR^3$:
 \begin{enumerate}
\item[(A)] For each $1>\alpha>0$, $\Sigma_j\setminus \cS$
converges in the $C^{\alpha}$-topology to the lamination $\cL
\setminus \cS$. \item[(B)]  $\sup_{B_{r}(x)\cap
\Sigma_j}|A|^2\to\infty$ as $j \to \infty$ for all $r>0$ and $x\in
\cS$.  (The curvatures blow up along $\cS$.)
\end{enumerate}
\end{theorem}

\subsection{Planar domains: the fine structure theorems}

  We will assume
here that the $\Sigma_i$'s are not disks  (recall that the case of
disks was dealt with in \cite{CM3}--\cite{CM6}). In particular, we
will assume that for each $i$, there exists some $y_i \in \RR^3$
and $s_i>0$ so that
\begin{equation}    \label{e:notulsc}
{\text{ some component of }} B_{s_i}(y_i) \cap \Sigma_i {\text{
 is not a disk.}}
\end{equation}
Moreover,  if the non-simply connected balls $B_{s_i}(y_i)$ ``run
off to infinity'' (i.e., if each connected component of
$B_{R_i'}(0) \cap \Sigma_i$ is a disk for some $R_i' \to \infty$),
then the results of \cite{CM3}--\cite{CM6} apply.  Therefore,
 after passing to a subsequence, we can assume that the
surfaces are uniformly not  disks, namely, that  there exists some
$R>0$ so that  \eqr{e:notulsc} holds with $s_i=R$ and $y_i=0$  for
all $i$.

In view of Theorem \ref{c:main} and the earlier results for disks,
it is natural to first analyze sequences that are ULSC, so where
$\cSt=\emptyset$, and second analyze sequences where $\cSu$ is
empty.  We will do this next.

\subsection{ULSC sequences}

Loosely speaking, our next result  shows that when the sequence is
ULSC (but not simply connected), a subsequence converges to a
foliation by parallel planes away from two lines $\cS_1$ and
$\cS_2$. The lines $\cS_1$ and $\cS_2$ are disjoint and orthogonal
to the leaves of the foliation and the two lines are precisely the
points where the curvature is blowing up. This is similar to the
case of disks, except that we get two singular curves for
non-disks as opposed to just one singular curve for disks.

\begin{theorem}[ULSC compactness, \cite{CM7}]  \label{t:t5.1}
 Let  a sequence $\Sigma_i$, limit lamination $\cL$, and singular
 set $\cS$
be  as in Theorem \ref{t:tab}. Suppose that each $\Sigma_i$
satisfies
 \eqr{e:notulsc}  with
$s_i=R > 1$ and $y_i=0$. If every $\Sigma_i$ is ULSC and
\begin{equation}
\sup_{B_1\cap \Sigma_i}|A|^2\to \infty \, ,
 \end{equation}
 then the limit lamination $\cL$ is the foliation $\cF = \{ x_3 = t\}_t$
  and the singular set $\cS$ is the union of two  disjoint
lines $\cS_1$ and $\cS_2$
 such that:
 \begin{enumerate}
\item[($C_{ulsc}$)]
Away from $\cS_1 \cup \cS_2$, each $\Sigma_j$ consists of exactly
two multi-valued graphs spiraling together.
 Near $\cS_1$ and $\cS_2$, the pair of multi-valued graphs form double
spiral staircases with opposite orientations at $\cS_1$ and
$\cS_2$. Thus, circling only $\cS_1$ or only $\cS_2$  results in
going either up or down, while a path circling both $\cS_1$ and
$\cS_2$ closes up. \item[($D_{ulsc}$)] $\cS_1$ and $\cS_2$ are
orthogonal to the leaves of the foliation.
\end{enumerate}
\end{theorem}

\begin{remark}
See Appendix \ref{s:s1} for the definition of a multi-valued
graph.  Roughly speaking a multi-valued graph is locally a graph
over a subset of a plane, but fails to be a global graph since the
projection to the plane is not one-to-one.
\end{remark}

\subsection{Sequences that are not ULSC}

When the sequence is no longer ULSC, one can get other types of
curvature blow-up by considering the family   of embedded minimal
planar domains known as the Riemann examples.
 Recall that, modulo translations and rotations, this is a
two-parameter family of periodic minimal surfaces, where the
parameters can be thought of as the size of the necks and the
angle  from one fundamental domain to the next.

With these examples in mind, we are now ready to state our second
main structure theorem describing the case where $\cSu$ is empty.

\begin{theorem}[\cite{CM7}] \label{t:t5.2a}
 Let  a sequence $\Sigma_i$, limit lamination $\cL$, and singular
 set $\cS$
be  as in Theorem \ref{t:tab}.  If $\cSu = \emptyset$ and
\begin{equation}
\sup_{B_1\cap \Sigma_i}|A|^2\to \infty \, ,
 \end{equation}
   then $\cS=\cSt$ by \eqr{e:} and
 \begin{enumerate}
\item[($C_{neck}$)]
Each point $y$ in $\cS$ comes with a sequence of
{\underline{graphs}} in $\Sigma_j$ that converge  to the plane $\{
x_3 = x_3 (y) \}$. The convergence is in the $C^{\infty}$ topology
away from the point $y$ and possibly also one other point  in $\{
x_3 = x_3 (y) \} \cap \cS$.  If the convergence is   away from one
point, then these graphs are defined over annuli; if the
convergence is away from two points, then the graphs are defined
over disks with two subdisks removed.
\end{enumerate}
\end{theorem}

\subsection{An overview of the proofs: The ULSC case}

  A key point
will be  that the results of \cite{CM3}--\cite{CM6} for disks will
give a sequence of multi-valued graphs in the $\Sigma_j$'s
 near each point
$x \in \cSu$.  Moreover, these multi-valued graphs close up in the
limit to give a leaf of $\cL'$ which extends smoothly across $x$.
Such a leaf is said to be {\emph{collapsed}}; in a neighborhood of
$x$, the leaf can be thought of as a limit of double-valued graphs
where the upper sheet collapses onto the lower. We   show that
every collapsed leaf is stable, has at most two points of $\cSu$
in its closure, and these points are removable singularities.
These results on collapsed leaves are applied first in the USLC
case
 and then again to get the structure of the ULSC
regions of the limit in general, i.e., (C2) and (D) in Theorem
\ref{t:t5.2}.

   Roughly speaking, there are
two main steps to the proof of Theorem \ref{t:t5.1}:
\begin{enumerate}
\item
Show that each collapsed leaf is in fact a plane punctured at two
points of $\cS$ and,  moreover, the sequence has the structure of
a double spiral staircase near both of these points, with opposite
orientations at the two points.
\item
Show that leaves which are nearby a collapsed leaf of $\cL'$ are
also planes punctured at two points of $\cS$.  (We call this
``properness''.)
\end{enumerate}

\subsection{An overview of the proofs: The general structure}

  Theorem \ref{t:t5.2a}, as well as Theorem \ref{c:main},
are proven by first analyzing sequences of minimal surfaces
without any assumptions on the sets $\cSu$ and $\cSt$.  The
precise statement of this general theorem is given in Appendix
\ref{s:genc}.  We will give an overview of the theorem next.

In this general case, we show that a subsequence
  converges to a lamination $\cL'$ divided into regions where
  Theorem \ref{t:t5.1} holds and  regions
  where Theorem \ref{t:t5.2a} holds.  This
convergence is in $C^{1,1}$ topology away from the singular set
$\cS$ where the curvature blows up.  Moreover, each point of $\cS$
comes with a plane and these planes are essentially contained in
$\cL'$. The set of heights of the planes is a closed subset $\cI
\subset \RR$ but may not be all of $\RR$ as it was in Theorem
\ref{t:t5.1} and may not even be connected. The behavior of the
sequence is different at the two types of singular points in $\cS$
- the set $\cSt$ of ``catenoid points'' and the set $\cSu$ of ULSC
singular points. We will see that $\cSu$ consists of a union of
Lipschitz curves transverse to the lamination $\cL$. This
structure of $\cSu$ implies that the set of heights in $\cI$ which
intersect $\cSu$ is a union of intervals; thus this part of the
lamination is foliated. In contrast, we will not get any structure
of the set of ``catenoid points'' $\cSt$. Given a point $y$ in
$\cSt$, we will get a sequence of graphs in $\Sigma_j$ converging
to a plane through $y$. This convergence will be in the smooth
topology away from either one or two singular points, one of which
is $y$. Moreover, this limit plane through $y$ will be a leaf of
the lamination $\cL$.

  The key steps for
proving the general structure theorem are the following:
\begin{enumerate}
\item Finding a {\underline{stable}} plane through each point of
$\cSt$.  This plane will be a limit of a sequence of stable
graphical annuli that lie in the complement of the surfaces. \item
Finding graphs in $\Sigma_j$ that converge to a plane through each
point of $\cSt$.  To do this, we look in regions between
consecutive necks and show that in any such region the surfaces
are ULSC.  The one-sided curvature estimate will then allow us to
show that these regions are graphical. \item Using (1) and (2) we
then analyze the ULSC regions of a limit. That is, we show that if
the closure of a leaf in $\cL'$ intersects $\cSu$, then it has a
neighborhood that is ULSC.  This will allow us to use the argument
for the proof of Theorem \ref{t:t5.1} to get the same structure
for such a neighborhood as we did in case where the entire
surfaces where ULSC.
\end{enumerate}

  The main point left in Theorem \ref{t:t5.2a}, which is not
  included in this general compactness theorem,
 is to prove that {\underline{every}} leaf of
  the lamination
  $\cL$ in Theorem \ref{t:t5.2a} is a plane.  In contrast,
the general compactness theorem gives a plane through each point
of $\cSt$, but does not claim that the leaves of $\cL'$ are
planar.

 Finally, since
the no mixing theorem implies that Theorem \ref{t:t5.1} and
Theorem
 \ref{t:t5.2a} cover all cases,
   Theorem \ref{t:tab} will be a corollary of
 these two theorems.

\section{The structure of embedded minimal annuli}

 We turn next to a local structure theorem for embedded minimal annuli that, roughly speaking, shows that
they must look like catenoids.  Namely, the main theorem of
\cite{CM9}  proves
 that any embedded minimal annulus in a ball
(with boundary in the boundary of the ball and) with a small neck
 can be decomposed by a simple closed geodesic into two graphical
sub--annuli.  Moreover, there is a sharp bound for the length of
this closed geodesic in terms of the separation (or height)
between the graphical sub--annuli.   This  serves to illustrate
the ``pair of pants'' decomposition from \cite{CM7} in the special
case where the embedded minimal planar domain is an annulus.

The precise statement of this decomposition for annuli is:

\begin{theorem}[Main Theorem, \cite{CM9}] \label{t:nitsche}
There exist $\epsilon>0$, $C_1 ,\, C_2,\, C_3>1$ so:   If
$\Sigma\subset B_{R}\subset \RR^3$ is an embedded minimal annulus
with $\partial \Sigma\subset \partial B_{R}$ and
$\pi_1(B_{\epsilon R}\cap \Sigma)\ne 0$, then there is a simple
closed geodesic $\gamma \subset \Sigma$ of length $\ell$ so that:

\begin{itemize}
\item
 The curve $\gamma$ splits the connected component of $B_{R/C_1}\cap
\Sigma$ containing it into
  annuli
$\Sigma^{+}$ and $ \Sigma^{-}$, each with $\int |A|^2 \leq 5 \,
\pi$.
 \item
 Each of $\Sigma^{\pm}
\setminus \cT_{C_2 \, \ell}(\gamma)$ is a graph with gradient
$\leq 1$. \item  $\ell \log (R/\ell) \leq C_3\,h$ where the
separation $h$ is given by \begin{equation}
 h=\min_{x_{\pm} \in\partial B_{
R/C_1 }\cap \Sigma^{\pm}}|x_+ -x_- | \, .
\end{equation}
\end{itemize}
\end{theorem}

Here  $\cT_{s}(S) \subset \Sigma$ denotes the intrinsic
$s$--tubular neighborhood of a subset $S\subset \Sigma$.

\subsection{A sketch of the proof}

We will next give a brief sketch of the proof of the decomposition
theorem, Theorem \ref{t:nitsche}.  The starting point is to use
the hypothesis $\pi_1(B_{\epsilon R}\cap \Sigma)\ne 0$ and a
barrier argument to find a stable graph $\Gamma_0$ that is defined
over an annulus and disjoint from $\Sigma$.  The stable graph
$\Gamma_0$ will allow us to divide $\Sigma$ into two pieces, one
on each   side of $\Gamma_0$.   To do this, we first fix a simple
closed $\tilde \gamma \subset B_{\epsilon R}\cap \Sigma$ that
separates the two boundary components of $\Sigma$. The curve
$\tilde \gamma$ is contained in a small extrinsic ball, but there
is no a priori reason why it must be short.{\footnote{However, the
chord arc bounds in the later paper \cite{CM16} could now be used
to bound the length.}} A barrier argument using a result of Meeks
and Yau then gives a stable embedded minimal annulus $\Gamma$ that
separates the two boundary components of $\Sigma$ and where
$\tilde \gamma$ is one component of the boundary $\partial \Gamma$
and the other component is in $\partial B_R$.  Finally, Theorem
$0.3$ of \cite{CM5} then implies that $\Gamma$ contains the
desired graph $\Gamma_0$; this should be compared with the
well-known result of D. Fischer-Colbrie in the complete case.

We will see next that each half of $\Sigma$, i.e., the part above
$\Gamma_0$ and the part below $\Gamma_0$, is itself a graph away
from the boundary of $\Gamma_0$.   This part of the argument
applies more generally to an ``annular end'' of a minimal surface.
We will prove that each half of $\Sigma$ contains a graph by
showing that it must contain large locally graphical pieces and
then using embeddedness to see that these pieces must be global
graphs (i.e., the projection down is one-to-one). This follows by
combining three facts:
\begin{enumerate}
\item[(1)]  The one-sided curvature estimate of
\cite{CM3}--\cite{CM6}
 gives a scale-invariant
curvature estimate for $\Sigma$'s in a narrow cone about the graph
$\Gamma_0$.  This requires that we know that each component of
$\Sigma$ in balls away from the origin is a disk; this can be seen
from the maximum principle. \item[(2)] Using (1), the gradient
estimate
 gives a  narrower cone about $\Gamma_0$ where
$\Sigma$ is locally graphical.  This is because (1) implies that
the surface is well-approximated by its tangent plane and, since
it cannot cross $\Gamma_0$, it must be almost parallel to
$\Gamma_0$. \item[(3)] As long as $\epsilon$ is small enough, each
half of $\Sigma$ must intersect any narrow cone about $\Gamma_0$.
This was actually proven in lemma $3.3$ of \cite{CM8} that gave
the existence of low points in a {\underline{connected}} minimal
surface contained on one side of a plane and with interior
boundary close to this plane.{\footnote{The argument for this was
by contradiction. Namely, if there were no low points, then we
would get a contradiction from the strong maximum principle by
first sliding a catenoid up under the surface and then sliding the
catenoid horizontally away, eventually separating two boundary
components of the surface. Here the strong maximum principle is
used to keep the sliding catenoids and the surface disjoint.  See,
for instance, corollary $1.18$ in \cite{CM1} for a precise
statement of the strong maximum principle.}} \end{enumerate}
 Step (3) allows us to find very flat regions in $\Sigma$ near
 $\Gamma_0$, we can then repeatedly apply the Harnack inequality to build
this out into large locally graphical regions that stay inside the
narrow cone about $\Gamma_0$. These locally graphical regions
piece together to give a graph over an annulus; the other
possibility would be to form a multi-valued graph, but this is
impossible since such a multi-valued graph would be forced to
spiral infinitely (since it cannot cross itself and also cannot
cross the stable graph $\Gamma_0$).

Finally, the last step of the proof is to use a blow up argument
to get the precise bounds on the length of the curve $\gamma$.

\subsection{Complete properly embedded minimal annuli}

The decomposition of properly embedded minimal annuli given by
Theorem \ref{t:nitsche} can be viewed as a local version of
well-known global results of P. Collin, \cite{Co}, and Colding and
the author, \cite{CM14}, on annular ends.

To explain these global results, recall that $\Sigma$ is said to
have {\emph{finite topology}} if it is homeomorphic to a closed
Riemann surface with a finite number of punctures; the genus of
$\Sigma$ is then the genus of this Riemann surface and the number
of punctures is the number of ends. It follows that a neighborhood
of each puncture corresponds to a properly embedded annular end of
$\Sigma$. Perhaps surprisingly at first, the more restrictive case
is when $\Sigma$ has more than one end. The reason for this is
that a barrier argument gives a stable minimal surface between any
pair of ends. Such a stable surface is then asymptotic to a plane
(or catenoid), essentially forcing each end to live in a
half--space. Using this restriction, P. Collin proved:

\begin{theorem}[Main theorem,  \cite{Co}]      \label{t:co}
Each end of a complete properly embedded minimal surface with
finite topology and at least two ends is asymptotic to a plane or
catenoid.
\end{theorem}

In particular, such a $\Sigma$ has finite total curvature and,
outside some compact set, $\Sigma$ is given by a finite collection
of disjoint graphs over a common plane.

As mentioned above, Collin proved Theorem \ref{t:co} by showing
that an embedded annular end that lives in a half--space must have
finite total curvature.  \cite{CM14} used the one--sided curvature
estimate to strengthen this from a half--space to a strictly
larger  cone, and in the process give a very different proof of
Collin's theorem.

\begin{theorem}[Main theorem,    \cite{CM14}]      \label{t:coCM}
There exists $\epsilon > 0$ so that any complete properly embedded
minimal annular end contained in the cone
\begin{equation}
    \{ x_3 \geq - \epsilon \, (x_1^2 + x_2^2 + x_3^2)^{1/2}  \}
\end{equation}
 is asymptotic to a plane or catenoid.
\end{theorem}

\section{Properness and removable singularities for minimal
laminations}    \label{s:proper}

The compactness theorems of \cite{CM3}--\cite{CM7} assume that the
surface $\Sigma_i$ has boundary $\partial \Sigma_i$ in the
boundary $\partial B_{R_i}$ of an expanding sequence of balls
where $R_i$ goes to infinity.  We call this the {\emph{global}}
case, in contrast to the {\emph{local}} case where the boundaries
are in the boundary of a fixed ball $\partial B_R$.{\footnote{One
can also consider the more restrictive {\emph{complete}} case
where $\Sigma_i$ is complete without boundary.}}

This distinction between the local and global cases explains why
 the global compactness theorem for sequences of
 disks does not imply the compactness theorem for ULSC sequences.
 Namely, even though the
 ULSC sequence consists {\underline{locally}} of disks, the compactness
result for disks
 was in the {\underline{global}} case where the radii go to
 infinity and hence does not apply.

In order to focus the discussion, we will explain the differences
between the global and local cases for disks.  The assumption that
$R_i \to \infty$ is used in the compactness theorem for disks in
two ways:
\begin{enumerate}
\item[(1)] We show that the limit lamination contains a stable
leaf through each singular point.  Since $R_i \to \infty$, this
stable leaf is complete and, hence,  a plane by the Bernstein
theorem of D. Fischer-Colbrie and R. Schoen and M. Do Carmo and C.
Peng. \item[(2)]  We show next that the leaves nearby a singular
point must also be planes.  It follows that the singular set
cannot stop and all of $\RR^3$ is foliated by planes in the limit.
We call this {\emph{properness}}.
\end{enumerate}

The use of $R_i \to \infty$ in (1) is not really essential. The
leaf would no longer have to be flat in the local case, but it
would satisfy uniform estimates  by R. Schoen's curvature estimate
for stable surfaces, \cite{Sc1} (cf. \cite{CM2}).

In contrast, it turns out that the use of $R_i \to \infty$ in (2)
is essential.  Namely, in  \cite{CM15}, we constructed
  a sequence of embedded minimal disks $\Sigma_i$ in
the unit ball $B_1$ with $\partial \Sigma_i \subset
\partial B_1$ where
 the curvatures blow up only at $0$
and
\begin{equation}
    \Sigma_i  \setminus \{ x_3 = 0 \}
\end{equation}
 converges to two
  embedded minimal disks
 \begin{align}
    \Sigma^- &\subset \{ x_3 < 0 \}  \\
 \Sigma^+ &\subset \{ x_3 > 0 \} \, ,
 \end{align}
  each of which spirals
 into $\{ x_3 =
0 \}$ and thus is not proper. Thus, in the example from
\cite{CM15}, $0$ is the first, last, and only point in $\cSu$ and
the limit lamination consists of three leaves: $\Sigma^+$,
$\Sigma^-$, and the punctured unit disk $B_1 \cap \{ x_3 = 0 \}
\setminus \{ 0 \}$.  This  lamination of
 $B_1 \setminus \{ 0\}$
 cannot be extended smoothly to a lamination of $B_1$; that is to say, $0$ is not a removable
 singularity.  This should be contrasted with the global case
 where every singular point is a removable singularity
 for the limit foliation by parallel planes.
B. Dean has constructed similar examples where the singular set is
an arbitrary finite set of points in the vertical axis; see
\cite{De}.

\subsection{A sketch of the proof of properness for disks}
To explain the proof of properness in the global case for disks,
we first need to see what could go wrong.  Suppose therefore that
the origin $0$ is a singular point and $\{ x_3 = 0 \}$ is the
corresponding limit plane.  It follows from the one-sided
curvature estimate that the intersection of each $ \Sigma_j$ with
a low cone about $\{ x_3 = 0 \}$ consists of two multi-valued
graphs for $j$ large (the fact that there are exactly two is
established in proposition II.1.3 in \cite{CM6}). There are now
two possibilities:
\begin{enumerate}
\item[(P)] The multi-valued graphs in the complement of the cone
close up in the limit to a foliation. \item[(N-P)] These
multi-valued graphs converge to a collection of graphs  and
{\emph{at least one}} multi-valued graph that spirals infinitely
on one side of $\{ x_3 = 0 \}$.
\end{enumerate}
As we saw above, the second case (N-P) can occur in the local
case.  We will explain why it cannot happen in the global case.

Suppose therefore that (N-P) holds and the limit contains a
multi-valued graph that spirals infinitely down to the plane $\{
x_3 = 0 \}$; it is the graph of a multi-valued function $u (\rho ,
\theta)$ defined for all $\rho \geq \ee$ and {\underline{all}}
$\theta
> 0$.  The separation $w(\rho,\theta)$ between consecutive sheets is by
definition
\begin{equation}
w(\rho,\theta)=u(\rho,\theta+2\pi)-u(\rho,\theta)\, .
\end{equation}
Since the limit is embedded and spirals downward, we must have  $w
< 0$.  We will actually work with the conformally changed
functions $\tilde u(x+iy)=u(\ee^x,y)$ and $\tilde
w(x+iy)=w(\ee^x,y)$ that are defined on the quadrant $\{ x
> 1 , \, y> 0 \}$.  The key point in the proof of properness is to
show that: \begin{enumerate} \item[(Key)] The vertical flux across
$\{ x = 1 \}$ is negative infinity. \end{enumerate}

{\underline{Why (Key) leads to a contradiction}}: To see why (Key)
leads to a contradiction, we need to recall more about the limit
in case (N-P).  Namely, we showed in \cite{CM6} that there must be
two multi-valued graphs spiralling together just as occurs for the
helicoid.  The same argument applies to both multi-valued graphs,
so both have unbounded negative flux across $\{ x = 1 \}$, i.e.,
over the circle of radius $\ee$ in the plane.  Moreover, we also
showed in \cite{CM6} that these two halves can be joined together
by  curves in the embedded minimal disk with a uniform bound on
the length of the curves.  For example, the helicoid contains two
infinite valued graphs and these can be connected by horizontal
lines.  In any case, this leads to a flux contradiction: Stokes'
theorem implies that the sum of the fluxes across compact
subcurves over the circle of radius $\ee$ must be bounded by the
length of the connecting curves.  However, the length of the
connecting curves is uniformly bounded and the fluxes across the
other curves both go to negative infinity.

{\underline{The idea of the proof of (Key)}}: In order to keep
things simple, we will pretend   that $u$ is a harmonic function;
this  serves to illustrate the main ideas. Since the separation
$w$ is locally the difference of two harmonic functions,   $w$ is
also harmonic; hence, the conformally changed functions $\tilde u$
and $\tilde w$ are harmonic on the quadrant $\{ x > 1 , \, y > 0
\}$.  Note that $\tilde u$ is positive and $\tilde w$ is negative.

The property (Key) is now roughly equivalent to showing that
\begin{equation}    \label{e:key}
    \int_0^{\infty} \frac{\partial \tilde u}{\partial x} (1,y) \,
    dy = + \infty \, .
\end{equation}
It may be helpful to consider an example; the function $\tilde u =
\pi / 2 - \arctan (y/x)$ is positive, harmonic, and its
multi-valued graph is an embedded infinite spiral that accumulated
to the plane $\{ x_3 = 0 \}$.  Furthermore, is is easy to verify
\eqr{e:key} in this case:
\begin{equation}    \label{e:keyex}
    \int_0^{\infty} \frac{\partial \tilde u}{\partial x} (1,y) \,
    dy = \int_0^{\infty} \frac{y}{1 + y^2}   \,
    dy =
    + \infty \, .
\end{equation}

To prove \eqr{e:key} for a general function $\tilde u$, first
observe that Stokes' theorem gives
\begin{equation}    \label{e:key}
    \int_1^{R} \frac{\partial \tilde u}{\partial x} (1,y) \,
    dy  + \int_1^{ R } \frac{\partial \tilde u}{\partial y} (x,1) \,
    dx = \int_{ \{ x^2 + y^2 = R^2 + 1 , \, x> 1, \, y> 1\} }
   \frac{ \partial  u}{\partial r} \, .
\end{equation}
In \cite{CM9}, we used the lower bound $\tilde u \geq 0$ to prove
that
\begin{equation}    \label{e:sphf}
\int_{ \{ x^2 + y^2 = R^2 + 1 , \, x> 1, \, y> 1\} }
   \frac{ \partial  u}{\partial r}
\end{equation}
is essentially non-negative.  The reason for this is that
\eqr{e:sphf} measures the logarithmic rate of growth of the
average of $\tilde u$ on the semi-circle; if this was negative,
the function $\tilde u$ would eventually have to become negative.

We then proved the claim \eqr{e:key}  in \cite{CM9} by using a
sharp estimate for the decay of $\tilde w$ to show that
\begin{equation}
\int_1^{ \infty } \frac{\partial \tilde u}{\partial y} (x,1) \,
    dx = - \infty \, .
\end{equation}
To explain this, observe that $\tilde w (x,1)$ is nothing more
than  $\tilde u (x,1+ 2\pi) - \tilde u (x,1)$ and, hence, can be
written as
\begin{equation}
    \tilde w (x,1) = \int_1^{1+ 2\pi} \frac{\partial \tilde u}{\partial
    y}(x,y) \, dy \, .
\end{equation}
In particular, $\tilde w (x,1)  \approx 2\pi \, \frac{\partial
\tilde u}{\partial y} (x,1)$.
 We proved in \cite{CM9} that
the fastest possible decay for $|\tilde w(x,1)|$ is $c_1/x$ and,
consequently, we get that
\begin{equation}
\int_1^{ \infty }    \tilde w (x,1) \,
    dx = - \infty \, .
\end{equation}
This completes the sketch of the proof.  The actual argument in
\cite{CM9} is somewhat more complicated, but similar in flavor.

\section{The uniqueness of the helicoid}

The helicoid and plane are the only classical examples of properly
embedded complete  minimal disks in $\RR^3$.
  It turns out that there is a good reason for the scarcity of examples.  Namely,
  using
the compactness theorem and one-sided curvature estimate of
\cite{CM3}--\cite{CM6}, W. Meeks and H. Rosenberg proved the
uniqueness of the helicoid:

\begin{theorem}     [Main theorem, \cite{MeRo1}]     \label{t:mero}
The plane and helicoid are the only complete properly embedded
simply--connected minimal surfaces in $\RR^3$.
\end{theorem}

This uniqueness has many applications, including additional
regularity of the singular set $\cS$.  To set this us, recall that
if we take a sequence of rescalings of the helicoid, then the
singular set $\cS$ for the convergence is the vertical axis
perpendicular to the leaves of the foliation. In \cite{Me1}, W.
Meeks used this fact together with the uniqueness of the helicoid
to prove that the singular set $\cS$ in Theorem \ref{t:t0.1} is
always a straight line perpendicular to the foliation.

There is an analog of Theorem \ref{t:mero} in the higher genus
case.  Namely, any properly embedded minimal surface with finite
(non-zero) genus and one end must be asymptotic to a helicoid.
Until recently, it was not known whether any such surface exists;
however, the construction of the genus one helicoid in
\cite{HoWeWo1} suggests that there may be a substantial theory of
these.

\begin{remark}
It follows from \cite{CM16} that any complete embedded
 minimal surface with finite topology in $\RR^3$ is automatically
properly embedded.  In particular, the hypothesis of  properness
can be removed from Theorem \ref{t:mero}.
\end{remark}

\subsection{A sketch of the proof} We will give a brief overview
of the proof by Meeks and Rosenberg for the uniqueness of the
helicoid; we refer to the original paper \cite{MeRo1} for the
details.

The first main step in the proof of Theorem \ref{t:mero} is to
analyze the asymptotic structure of a non-flat embedded minimal
disk $\Sigma$, showing that it looks roughly like a helicoid. This
is done in \cite{MeRo1} by analyzing sequences of rescalings of
$\Sigma$.  This  rescaling argument yields a sequence of embedded
minimal disks which does not converge in the classical sense
(there are no local area bounds). However, the lamination theorem
of \cite{CM3}--\cite{CM6} gives that a subsequence converges to a
foliation by parallel planes away from a Lipschitz curve
transverse to these planes.  Moreover, the lamination theorem also
gives that the intersection of $\Sigma$ with a cone consists of
two asymptotically flat multi-valued graphs. In particular,   this
foliation is unique, i.e., does not depend on the choice of
subsequence.  After possibly rotating $\RR^3$, we can assume that
the limit foliation is by horizontal planes (i.e., level sets of
$x_3$).

The second main step is to show that the height function $x_3$
together with its harmonic conjugate (which we will denote by
$x_3^{\ast}$) give global isothermal coordinates on $\Sigma$. This
step is crucial since it reduces the problem to analyzing the
potential Weierstrass data on the plane.  There are two key
components to getting these global coordinates, both of
independent interest. First, one must show that $\nabla x_3$ does
not vanish on $\Sigma$ -- i.e., that the Gauss map misses the
north and south poles.    Second, one must show that the map $(x_3
, x_3^{\ast})$ is proper -- i.e., that $x_3^{\ast}$ goes to
infinity as we go out horizontally.  Both of these steps strongly
use the asymptotic structure established in the first step.

The third main step is to analyze the Weierstrass data  in the
conformal coordinates $(x_3 , x_3^{\ast})$.   In these
coordinates, the only unknown is a  meromorphic function $g$ that
is the stereographic projection of the Gauss map of $\Sigma$.
 Since the Gauss map was already shown to miss the
north and south poles, the function $g$ can be written as $g=
\ee^{f}$ for an entire holomorphic function $f$. Meeks and
Rosenberg then use a Picard type argument   to show that $f$ must
be linear. The key in this argument is to analyze  the inverse
images of horizontal circles in $\SS^2$ under the Gauss map $\nn$,
using rescaling arguments and the compactness theory of
\cite{CM3}--\cite{CM6} to control the number of components.
Finally, every linear function $f$ gives rise to a surface in the
associate surface family of the helicoid, but the actual helicoid
is the only one of these that is embedded.

\section{Quasiperiodicity of properly embedded minimal planar domains}

We turn next to recent results of W. Meeks, J. Perez, and A. Ros
on the structure of  {\emph{complete}} properly embedded minimal
planar domains  with infinitely many ends in $\RR^3$. They have
obtained many important results on these surfaces in the series of
papers \cite{MePRs1}--\cite{MePRs3}; we have chosen to focus on
the main result of \cite{MePRs1}.

Many of these structure results are motivated by the two-parameter
family of minimal planar domains known as the Riemann examples;
see

\subsection{A few definitions}
To set this up, we first recall a few properties of
 a complete properly embedded minimal
planar domain $M \subset \RR^3$.

First, it follows from a barrier argument of Meeks and Yau that
one can find a stable embedded annulus between each pair of ends
of $M$; a result of Fischer-Colbrie then implies that this stable
surface has finite total curvature, so its ends are asymptotic to
planes or half-catenoids.  Since $M$ is embedded, these planes or
half-catenoids between its ends must all be parallel; this plane
is the {\emph{limit tangent plane at infinity}}.

Using these planes and half-catenoids in this way, Callahan,
Frohman, Hoffman, and Meeks showed that the ends of $M$ are
ordered by height over this limit tangent plane at infinity.
Moreover, a nice argument of Collin, Kusner, Meeks, and Rosenberg
in \cite{CKMR} shows that there are at most two {\emph{limit
ends}}; these can only be on the ``top'' or on the ``bottom''.

Finally, since the coordinate functions are harmonic on any
minimal surface, it follows from Stokes' theorem that the flux of
a coordinate function $x_i$ around a closed curve $\gamma$ depends
only on its homology class $[\gamma]$.  Recall that  the flux of
$x_i$ around $\gamma$ is
\begin{equation}
    \int_{\gamma} \frac {\partial x_i}{\partial n} \, ,
\end{equation}
where $\frac {\partial x_i}{\partial n}$ is the derivative of
$x_i$ in the conormal direction, i.e., the direction tangent to
$M$ but normal to $\gamma$. Using all three coordinates functions
at once gives the {\emph{flux}} map
\begin{equation}
    Flux : [\gamma] \to \RR^3 \, .
\end{equation}

\subsection{The curvature estimate of \cite{MePRs1}}
We can now define
 $\cM$ to be the space of properly embedded minimal planar domains in
$\RR^3$ with two limit ends, normalized so that every surface
$M\in\cM$ has horizontal limit tangent plane at infinity and the
vertical component of its flux equals one.  Here the horizontal
flux is the projection of the flux vector to the limit tangent
plane at infinity.

The main theorem of \cite{MePRs1} is:

\begin{theorem}[Theorem 5, \cite{MePRs1}]   \label{t:mpr1}
If a sequence $ M_i \in \cM$ has bounded horizontal flux, then the
Gaussian curvature of the sequence is uniformly bounded.
\end{theorem}

\begin{remark}
The main result in
 \cite{MePRs2} is that there must always be two limit ends; thus
 this hypothesis can be removed from Theorem \ref{t:mpr1}.
\end{remark}

   The first important application of the curvature
estimates is to describe the geometry of a properly embedded
minimal planar domain $M$ with  two limit ends:
\begin{quote}
 Such a surface $M$ has bounded curvature and is conformally a
compact Riemann surface punctured in a countable closed subset
with two limit points; the spacing between consecutive ends is
bounded from below in terms of the bound for the curvature; $M$ is
quasiperiodic in the sense that there exists a divergent sequence
$V(n)\in\Bbb R^3$ such that the translated surfaces $M+V(n)$
converge to a properly embedded minimal surface of genus zero, two
limit ends, a horizontal limit tangent plane at infinity and with
the same flux as $M$.
\end{quote}

Another particularly interesting consequence of Theorem
\ref{t:mpr1} is a solution of an old conjecture of Nitsche:

\begin{theorem}[Theorem 1,  \cite{MePRs1}]   \label{t:mpr2}
Any complete minimal surface which is a union of simple closed
curves in horizontal planes must be a catenoid.
\end{theorem}

\subsection{A heuristic argument for Theorem \ref{t:mpr1}}
 Theorem \ref{t:mpr1} is best illustrated by considering the
family of singly-periodic minimal surfaces known as the Riemann
examples. After normalizing so the vertical flux is one and
rotating so the horizontal flux points in the $x_1$-direction,
there is a $1$-parameter family parameterized by the length of the
horizontal flux $H$.  As $H \to 0$, this degenerates to a
catenoid; when $H \to \infty$, this degenerates to a helicoid.

With this in mind, there is  a simple heuristic argument for why
such an estimate  should hold.  Namely,  assume that $|A|^2(0) \to
\infty$ for a sequence $\Sigma_i$ and consider two possibilities:
\begin{enumerate}
\item
    The injectivity radius goes to zero.
\item
    Each $\Sigma_i$ is uniformly simply connected.
\end{enumerate}

In the first case, we would get short dividing curves in
$\Sigma_i$; integrating around these would then imply that the
flux was going to zero (violating the normalization).

In the second case, the results of Colding and the author in
\cite{CM7} give a limit plane through the origin which, by
 the uniqueness of the helicoid of Meeks and H. Rosenberg, is locally modelled by the helicoid
 for large $i$.  As discussed above, this corresponds (roughly, at least)   to a great
 deal of horizontal flux, violating the assumed bound on this.

The actual argument  is much more complicated but has at least a
similar flavor. One reason for the complications is that the
dichotomy (1) or (2) above is more subtle; (1) involves the
intrinsic distance while (2) uses the extrinsic distance. This
dichotomy can now be made rigorous using the proof of the
Calabi-Yau conjectures for embedded minimal surfaces with finite
topology in \cite{CM16}; however, the proof of Theorem
\ref{t:mpr1} came before \cite{CM16}, so intrinsic and extrinsic
distances were not yet known to be equivalent.

\appendix

\section{Multi-valued graphs}   \label{s:s1}

We have used two notions of  multi-valued graphs - namely, the one
used in \cite{CM3}--\cite{CM6} and a generalization.

In  \cite{CM3}--\cite{CM6}, we defined multi-valued graphs as
multi-sheeted covers of the punctured plane. To be precise, let
$D_r$ be the disk in the plane centered at the origin and of
radius $r$ and let $\cP$ be the universal cover of the punctured
plane $\CC\setminus \{0\}$ with global polar coordinates $(\rho,
\theta)$ so $\rho>0$ and $\theta\in \RR$. An $N$-valued graph of a
function $u$ on the annulus $D_s\setminus D_r$ is a single valued
graph over
\begin{equation}
     \{(\rho,\theta)\,|\,r\leq \rho\leq s\, ,\, |\theta|\leq
N\,\pi\} \, .
\end{equation}
Note that the helicoid is the union of two infinite--valued graphs
over the punctured plane together with the vertical axis.

 Locally, the above multi-valued graphs give the complete picture
 for a ULSC sequence.   However, the global picture can consist of
 several different multi-valued graphs glued together.  To
 allow for this, we are forced to consider multi-valued graphs
 defined over the universal cover of $\CC \setminus P$ where $P$
 is a discrete subset of the complex plane $\CC$.
 We will see that the bound on the
 genus implies that $P$ consists of at most two points.

\section{The lamination theorem and one-sided curvature estimate}

The first theorem that we recall shows that embedded minimal disks
are either graphs or are part of    double spiral staircases;
moreover, a sequence of such disks with curvature blowing up
converges to a foliation by parallel planes away from a singular
curve $\cS$.  This theorem is modelled on rescalings of the
helicoid and the precise statement is as follows (we state the
version for extrinsic balls; it was extended to intrinsic balls in
\cite{CM10}):

\begin{theorem} \label{t:t0.1}
(Theorem 0.1 in \cite{CM6}.)
 Let $\Sigma_i \subset B_{R_i}=B_{R_i}(0)\subset \RR^3$
be a sequence of embedded minimal {\underline{disks}} with
$\partial \Sigma_i\subset \partial B_{R_i}$ where $R_i\to \infty$.
If \begin{equation}
    \sup_{B_1\cap \Sigma_i}|A|^2\to \infty \, ,
    \end{equation}
 then there exists a subsequence, $\Sigma_j$, and
a Lipschitz curve $\cS:\RR\to \RR^3$ such that after a rotation of
$\RR^3$:
\begin{enumerate}
\item[\underline{1.}] $x_3(\cS(t))=t$.  (That is, $\cS$ is a graph
over the $x_3$-axis.)
\item[\underline{2.}]  Each $\Sigma_j$
consists of exactly two multi-valued graphs away from $\cS$ (which
spiral together).
\item[\underline{3.}] For each $1>\alpha>0$,
$\Sigma_j\setminus \cS$ converges in the $C^{\alpha}$-topology to
the foliation, $\cF=\{x_3=t\}_t$, of $\RR^3$.
\item[\underline{4.}] $\sup_{B_{r}(\cS (t))\cap
\Sigma_j}|A|^2\to\infty$ for all $r>0$, $t\in \RR$.  (The
curvatures blow up along $\cS$.)
\end{enumerate}
\end{theorem}

The second theorem that we need to recall asserts that every
embedded minimal disk lying above a plane, and coming close to the
plane near the origin, is a graph.   Precisely this is the {\it
intrinsic one-sided curvature estimate} which follows by combining
\cite{CM6} and \cite{CM10}:

\begin{theorem}  \label{t:t2}
There exists $\epsilon>0$, so that if
\begin{equation}
    \Sigma \subset \{x_3>0\} \subset \RR^3
\end{equation}
is an embedded minimal {\underline{disk}} with $\cB_{2R} (x)
\subset \Sigma \setminus
\partial \Sigma$ and $|x|<\epsilon\,R$,  then
\begin{equation}        \label{e:graph}
\sup_{ \cB_{R}(x) } |A_{\Sigma}|^2 \leq R^{-2} \, .
\end{equation}
\end{theorem}

Theorem \ref{t:t2} is in part used to prove the regularity of the
singular set where the curvature is blowing up.

\vskip2mm Note that the assumption in Theorem \ref{t:t0.1} that
the surfaces are disks is crucial and cannot even be replaced by
assuming that the sequence is ULSC. To see this, observe that  one
can choose a one-parameter family of Riemann examples which is
ULSC but where the singular set $\cS$ is given by a
{\underline{pair}} of vertical lines. Likewise, the assumption in
Theorem \ref{t:t2} that $\Sigma$ is simply connected is crucial as
can be seen from the example of a rescaled catenoid, see
\eqr{e:cat}.
 Under rescalings the catenoid converges (with
multiplicity two) to the flat plane. Thus a neighborhood of the
neck can be scaled arbitrarily close to a plane but the curvature
along the neck becomes unbounded as it gets closer to the plane.
Likewise, by considering the universal cover of the catenoid, one
sees that embedded, and not just immersed, is needed in Theorem
\ref{t:t2}.

\section{The precise statement of the general compactness theorem}
\label{s:genc}

 The precise statement of the compactness theorem for sequences
 that are neither necessarily ULSC nor with $\cSu = \emptyset$
 is the following:

\begin{theorem}[\cite{CM7}] \label{t:t5.2}
 Let $\Sigma_i \subset B_{R_i}=B_{R_i}(0)\subset \RR^3$
be a sequence of compact embedded minimal planar domains  with
$\partial \Sigma_i\subset \partial B_{R_i}$ where $R_i\to \infty$.
If
\begin{equation}
\sup_{B_1\cap \Sigma_i}|A|^2\to \infty \, ,
 \end{equation}
   then there is a subsequence $\Sigma_j$, a closed set $\cS$,
    and a lamination $\cL'$ of
$\RR^3 \setminus \cS$ so that:
\begin{enumerate}
\item[(A)] For each $1>\alpha>0$, $\Sigma_j\setminus \cS$
converges in the $C^{\alpha}$-topology to the lamination $\cL'$.
\item[(B)]  $\sup_{B_{r}(x)\cap \Sigma_j}|A|^2\to\infty$ as $j \to \infty$ for
all $r>0$ and $x\in \cS$.  (The curvatures blow up along $\cS$.)
\item[(C1)]
($C_{neck}$) from Theorem \ref{t:t5.2a} holds for each point
 $y$ in $\cSt$.
\item[(C2)]
($C_{ulsc}$) from Theorem \ref{t:t5.1} holds locally near $\cSu$.
 More precisely, each point $y$ in $\cSu$ comes with a sequence
of {\underline{multi-valued graphs}} in $\Sigma_j$ that converge
to the plane $\{ x_3 = x_3 (y) \}$. The convergence is in the
$C^{\infty}$ topology away from the point $y$ and possibly also
one other point   in $\{ x_3 = x_3 (y) \} \cap \cSu$.  These two
possibilities correspond to the two types of multi-valued graphs
defined in Section \ref{s:s1}.
\item[(D)]     The set $\cSu$ is a union of Lipschitz curves
transverse to the lamination.  The leaves intersecting $\cSu$ are
planes foliating an open subset of $\RR^3$ that does not intersect
$\cSt$. For the set $\cSt$, we make no claim about the structure.
\item[(P)]
Together (C1) and (C2) give a sequence of graphs or multi-valued
graphs converging to a plane through each point of $\cS$. If $P$
is one of these planes, then each leaf of $\cL'$ is either
disjoint from $P$ or is contained in $P$.
\end{enumerate}
\end{theorem}

 Note that Theorem \ref{t:t5.2}
 is a technical tool that is superseded
 by the stronger compactness theorems in the ULSC and non-ULSC
 cases, Theorem \ref{t:t5.1} and Theorem
 \ref{t:t5.2a}.  This is because   we will know by the no mixing
 theorem that either $\cSt =
\emptyset$ or $\cSu = \emptyset$, so that these cover all possible
cases.

\end{document}